\documentclass{birkjour}
 \newtheorem{thm}{Theorem}[section]

 \theoremstyle{definition}
 \newtheorem{defn}[thm]{Definition}
 \theoremstyle{remark}

 \numberwithin{equation}{section}

\begin{document}

%
%
%
%
%
%
%
%
%

\title[The Monge Theorem]
 {Does The Monge Theorem Apply To Some Non-Euclidean Geometries ? }

\author[Corresponding author]{Temel\ Ermi\c{s}*}

\address{%
Department of Mathematics and Computer Sciences,\\Eskisehir Osmangazi University, 26040 Eskisehir, Turkey}

\email{termis@ogu.edu.tr}

\thanks{* Corresponding author}
\author{\"{O}zcan Geli\c{s}gen}
\address{Department of Mathematics and Computer Sciences,\\Eskisehir Osmangazi University, 26040 Eskisehir, Turkey}
\email{gelisgen@ogu.edu.tr}


\date{January 1, 2004}

\begin{abstract}
In geometry, Monge's theorem states that for any three non-overlapping circles of distinct radii in the two dimensional analytical plane equipped with the Euclidean metric, none of which is completely inside one of the others, the intersection points of each of the three pairs of external tangent lines are collinear. So, it is clearly observed that Monge's theorem is an application of Desargues' theorem. Our main motivation in this study is to show whether Monge theorem is still valid even if the plane is equipped with the metrics alpha and  $L_{p}$.
\end{abstract}

\maketitle
\section{Introduction}
Euclidean distance is the most common use of distance. In most cases when
people said about distance, they will refer to Euclidean distance. The
Euclidean distance between two points is defined as the length of the
segment between two points. Although it is the most popular distance
function, it is not practical when we measure the distance which we actually
move in the real world ( we live on a spherical Earth rather than on a
Euclidean $3-$space \textbf{!} ). We must think of the distance as though a
car would drive in the urban geography where physical obstacles have to be
avoided. So, one had to travel through horizontal and vertical streets to
get from one location to another. To compensate disadvantage of the
Euclidean distance, the taxicab geometry was first introduced by K. Menger
and has developed by E. F. Krause using the taxicab metric $d_{T}$ of which
paths composed of the line segments parallel to coordinate axes (see Figure
1) (\cite{Krause}, \cite{Menger}). Its three-dimensional version has been
introduced in \cite{Akca}.
Later, researchers have wondered whether there are alternative distance
functions of which paths are different from path of Euclidean metric. For
example, G. Chen \cite{Chen} developed Chinese checker distance $d_{CC}$ in $%
\mathbb{R}^{2}\ $of which paths are similar to the movement made by Chinese
checker (see Figure 1).
Another example, S. Tian gave a family of metrics in \cite{Tian}, alpha
distance $d_{\alpha }$ for $\alpha \in \left[ 0,\pi /4\right] $, which
includes the taxicab and Chinese checker metrics as special cases. Then,
Kaya et al. have given the most general form of $d_{\alpha }$ on $n-$%
dimensional analytical space for $%
\alpha \in \left[ 0,\pi /2\right) $. When we examine the common features of
the metrics $d_{M}$, $d_{T}$,\ $d_{CC}\ $and$\ d_{\alpha }$, we see that
these metrics whose paths are parallel to at least one of the coordinate
axes. So, it is a logical question " Are there metric or metrics of which
paths are not parallel to the coordinate axes" ( Yes, there are \textbf{!}%
).\ The references \cite{Colakoglu-Alpha} and \cite{Polar}  
can be reviewed for the answer to this question. Meanwhile, the
alpha metric is defined as follows;
\begin{defn}
	Let $A=\left(x_{a},y_{a}\right) $ \ and \ $B=\left( x_{b},y_{b}\right) $ be two any points in $\mathbb{R}^{2}$ such that $\Delta _{AB}{\small =}\max \left\{
	\left\vert x_{a}-x_{b}\right\vert ,\left\vert y_{a}-y_{b}\right\vert
	\right\} $ and $\delta _{AB}{\small =}\min \left\{ \left\vert x_{a}{\small -}x_{b}\right\vert
	,\left\vert y_{a}{\small-}y_{b}\right\vert \right\} $. Then, the Euclidean metric can be given by 
	\begin{equation*}
		d_{E}\left( A,B\right) =\left( \Delta _{AB}^{2}+\delta _{AB}^{2}\right)
		^{1/2}\text{.}
	\end{equation*}%
	Also, the metrics $d_{M},\ d_{T}$,\ $d_{CC}\ $and$\ d_{\alpha }$ are defined
	as following Figure 1.
\end{defn}
\begin{figure}[h]
	\begin{center}
	\includegraphics[width=1.1\linewidth]{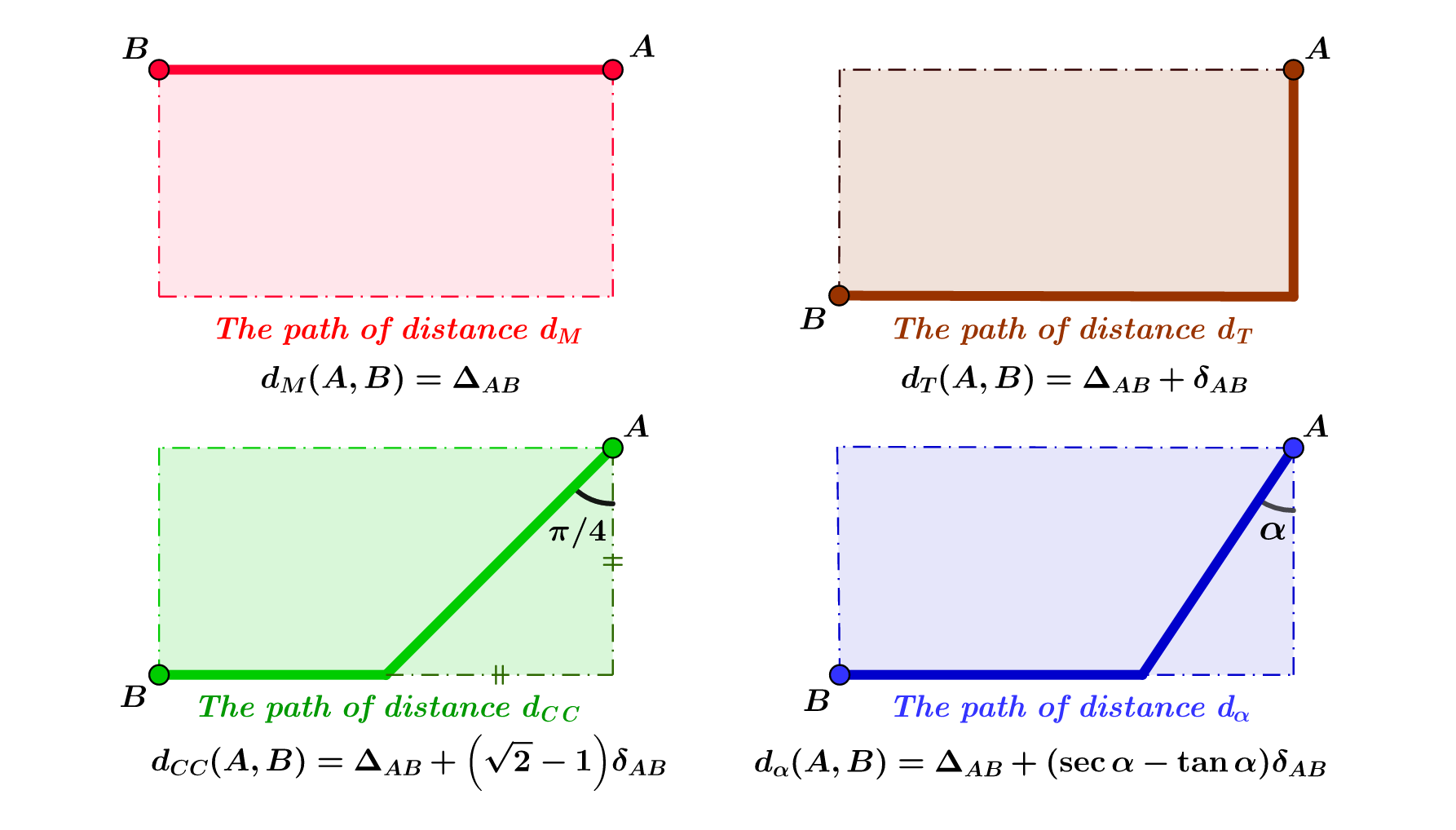}
	\caption{The paths of distances $d_{M},\ d_{T}$,\ $d_{CC}\ $and$\ d_{\alpha }$.}
	\label{fig:boat1}
	\end{center}
\end{figure}
 $\alpha -$(alpha) plane geometry, which includes the taxicab and Chinese
checker geo\-metry, is a Minkowski geometry. Minkowski geometry is a non
Euclidean geometry in a finite number of dimensions that is different from
elliptic and hyperbolic geometry (from Minkowskian geometry of space-time).
Here the linear structure is same as the Euclidean one but distance is not
uniform in all directions (see for details \cite{Thomas}, \cite{Lassak}, \cite{Swanepoel1}, \cite{Swanepoel2}, \cite{Mustafaev}, \cite{Thompson}). That is, $\alpha-$plane is almost the same as
Euclidean plane since the points are the same, the lines are the same, and
the angles are measured in the same way. Instead of the usual circle in
Euclidean plane geometry, unit ball is a certain symmetric closed (see
Figure 2). Since the $\alpha -$plane geometry has a different distance
function, it seems interesting to study the $\alpha -$analog of the topics
that include the concepts of distance in the Euclidean geometry. One of the
famous theorems that includes the concept of distance is the Monge Theorem \cite{Leopol}, \cite{Walker}.
Monge's theorem says that for any three nonintersecting circles in a plane,
none of which is equal radius, the intersection points of each of the three
pairs of external tangent lines are collinear. For any two circles in a
plane, an external tangent is a line that is tangent to both circles but
does not pass between them. There are two such external tangent lines for
any two circles. Each such pair has a unique intersection point in the
plane. In the case of two of the circles being of equal size, the two
external tangent lines are parallel. If the two external tangents are
considered to intersect at the point at infinity, then the other two
intersection points must be on a line passing through the same point at
infinity, so the line between them takes the same angle as the external
tangent.

We will show that Monge theorem is still valid even if the plane is equipped
with the metrics $\alpha$ and $L_{p}\ $in the next section. First, let us
recall the definition of the well-known $L_{p}-$metric.
\begin{defn}
	Let $X=\left( x_{1},\ldots
	,x_{n}\right) $ and $Y=\left( y_{1},\ldots ,y_{n}\right) $ be two vectors in
	the $n-$dimensional real vector space $\mathbb{R}^{n}$.\ The $l_{p}-$metric $%
	d_{l_{p}}$,\ $1\leq p\leq \infty $, is a norm metric on $\mathbb{R}^{n}$ (
	or on $\mathbb{C}^{n}$), is defined by 
	\begin{equation*}
		\left\Vert X-Y\right\Vert _{p}\text{,}
	\end{equation*}%
	where the $l_{p}-$norm\ $\left\Vert .\right\Vert _{p}$ is defined by 
	\begin{equation*}
		\left\Vert X\right\Vert _{p}=\left( \underset{i=1}{\overset{n}{\sum }}%
		\left\vert x_{i}\right\vert ^{p}\right) ^{\frac{1}{p}}\text{.}
	\end{equation*}%
	If $p\rightarrow \infty $, we obtain $\left\Vert X\right\Vert _{\infty }=%
	\underset{p\rightarrow \infty }{\lim }\left( \underset{i=1}{\overset{n}{%
			\sum }}\left\vert x_{i}\right\vert ^{p}\right) ^{\frac{1}{p}}=\underset{%
		1\leq i\leq n}{\max }\left\vert x_{i}\right\vert $.\ Also, the the pair $%
	\left( \mathbb{R}^{n},d_{l_{p}}\right) $ is called $l_{p}-$(metric) space.
\end{defn}
\section{The Unit Circles in $\mathbb{R}^{2}$ and Monge' s Theorem}
As is known to all, the circle is simply the set of points that are an equal
distance from a certain central point. That is, let $M$ be given point in
the plane, and $r$ be a positive real number. The set of points $\left\{
X\in \mathbb{R}^{2}:d\left( M,X\right) =r\right\} $ is called a circle, the
point $M$ is called center of the circle, and $r$ is called the length of
the radius or simply radius of the circle. But there's a big assumption in
definition of the circle: " distance ". While we usually use the Euclidean
distance $d_{E}$, there are other valid notions of distance we could have
used instead, which we have mentioned above.  
\begin{figure}[h]
	\begin{center}
		\includegraphics[width=0.8\linewidth]{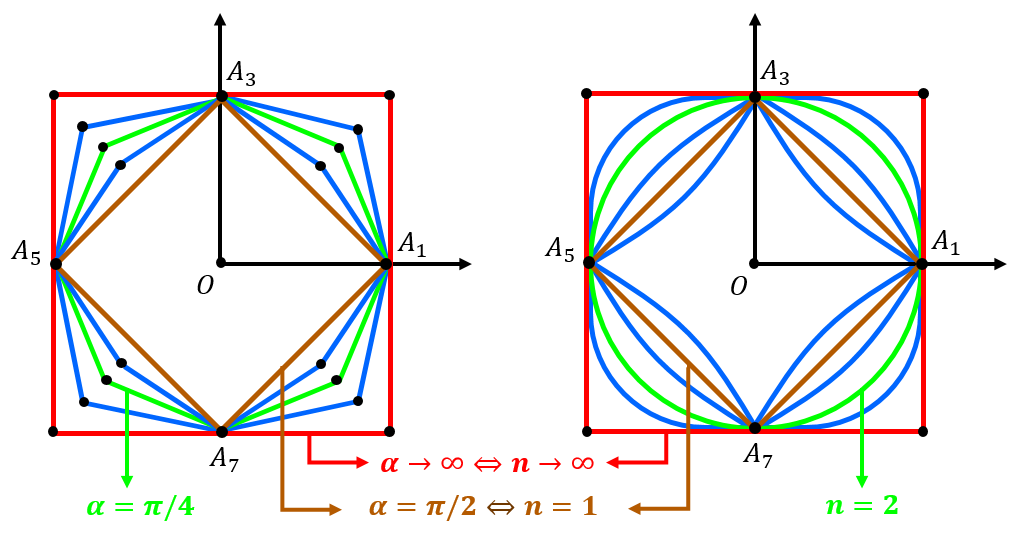}
		\caption{The unit circles in the $\alpha$ and $l_{p}-$plane.}
		\label{fig:boat1}
	\end{center}
\end{figure}
It is easy to see that if $\alpha \in \left( 0,\pi /2\right) $, then unit $\alpha-$circle is an octagon with corner points$\ A_{1}=\left( 1,0\right) $%
,\ $A_{2}=\left( 1/k,1/k\right) $,\linebreak $A_{3}{\small =}\left( 0,1\right) $,\ $%
A_{4}{\small =}\left(-1/k,1/k\right) $,\ $A_{5}{\small =}\left( -1,0\right) $,\ $A_{6}{\small =}\left(
-1/k,-1/k\right) $,\ $A_{7}=\left( 0,-1\right) $ and $A_{8}=\left(
1/k,-1/k\right) $, where$\ k=1+\sec \alpha -\tan \alpha $.$\ $Note that $%
A_{2}$ and $A_{6}$ are on the line $y=x$; $A_{4}$ and $A_{8}$ are on the
line $y=-x$. When $\alpha \rightarrow \pi /2$ and $\alpha =0$, the $\alpha-$circle is the circle with respect to the metrics $d_{M}$ and $d_{T}$,
respectively (see Figure 2). Also, $\alpha-$circle is the circle with
respect to the metrics $d_{CC}$ such that $\alpha =\pi /4$. Similarly, when $%
p\rightarrow \infty $,$\ p=1$\ and $p=2$, the circle in $l_{p}-$plane is the
circle with respect to the metrics $d_{M}$,$\ d_{T}$ and $d_{E}$,
respectively (see Figure 2). Now, when we take the circles $\alpha$ and $l_{p}$ instead of Euclidean circles, we will investigate whether the Monge theorem is still valid.
\begin{figure}[h]
	\begin{center}
		\includegraphics[width=0.84\linewidth]{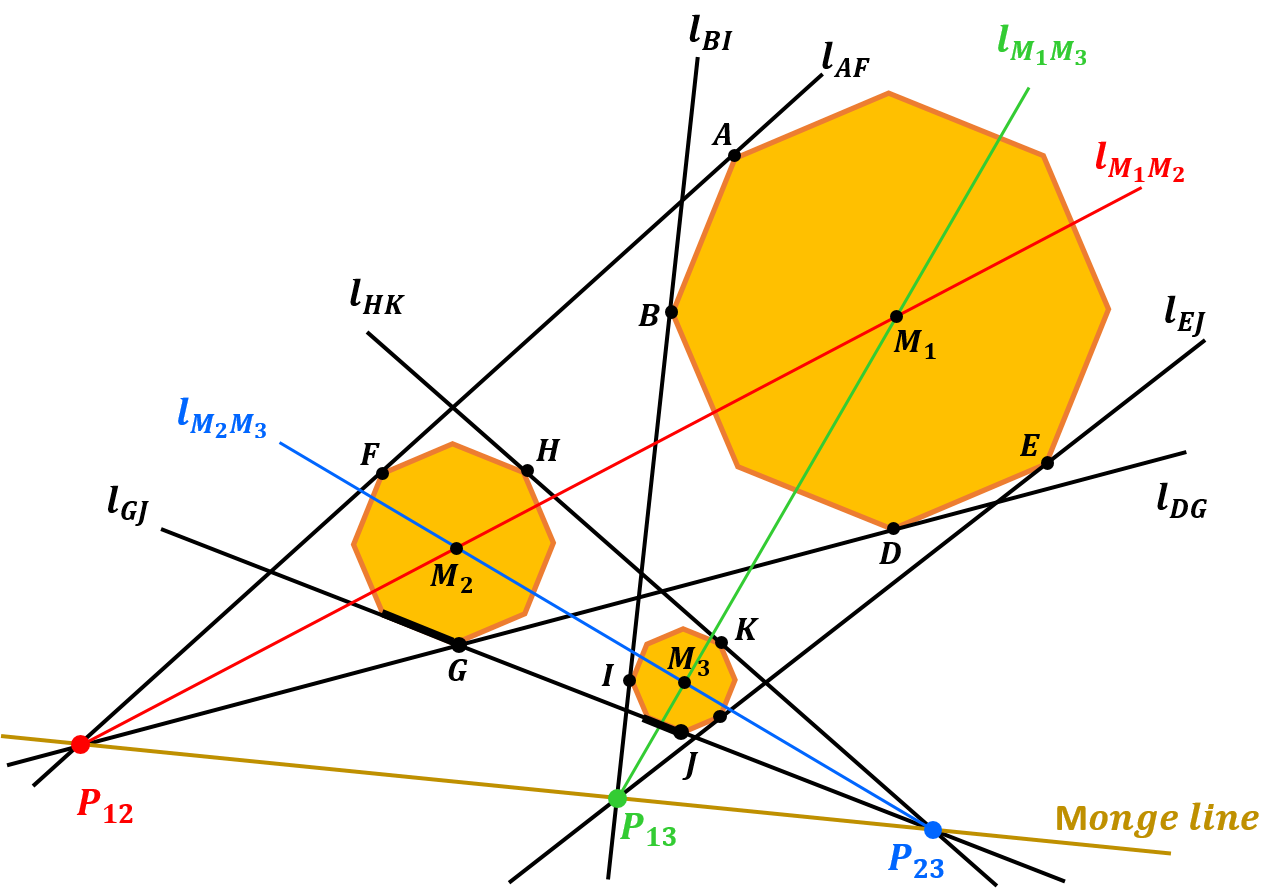}
		\caption{}
		\label{fig:boat1}
	\end{center}
\end{figure}
Consider Figure 3. Let $C_{i}$\ be the $\alpha -$circles (or $l_{p}-$circles ) with center \linebreak $%
M_{i}=\left( x_{i},y_{i}\right) $ and radius $r_{i}$ for $i=1,2,3$. Then, the
coordinates of points on the $\alpha -$circles (or $l_{p}-$circles )$\ C_{1}$%
, $C_{2}$ and $C_{3}$\ can be easily given with respect to the coordinates
of circles's centers and radii. Because we know the coordinates of points $%
A,F,D$ and $G$ on the circles, we can find the equations of the tangent
lines $l_{AF}$ and $l_{DG}$.\ Consequently, coordinates of the point $P_{12}$
of intersection of the tangent lines $l_{AF}$ and $l_{DG}\ $ be\ able to
calculated as\medskip\ $P_{12}=\left( \dfrac{r_{1}x_{2}-r_{2}x_{1}}{%
	r_{1}-r_{2}},\dfrac{r_{1}y_{2}-r_{2}y_{1}}{r_{1}-r_{2}}\right) $.\ Also, the
line $l_{M_{1}M_{2}}$ through the center points $M_{1}$ and $M_{2}$ passes
through the point $P_{12}$.\ Similarly, after calculating the equations of
tangent lines\medskip\ $l_{BI}$ and $l_{EJ}$ for the $\alpha -$circles\ (or $%
l_{p}-$circles )\ $C_{1}\ $and $C_{3}$, their intersection point of lines $%
l_{BI}$ and $l_{EJ}$ is found as $P_{13}=\left( \dfrac{r_{1}x_{3}-r_{3}x_{1}%
}{r_{1}-r_{3}},\dfrac{r_{1}y_{3}-r_{3}y_{1}}{r_{1}-r_{3}}\right) $. Finally,
coordinates\medskip\ of the point $P_{23}\ $be\ able to calculated as $%
P_{23}=\left( \dfrac{r_{2}x_{3}-r_{3}x_{2}}{r_{2}-r_{3}},\dfrac{%
	r_{2}y_{3}-r_{3}y_{2}}{r_{2}-r_{3}}\right) $. By simplifying  long and boring calculations, we get the equation of the line
passing through points $P_{12}$ and $P_{13}$ as 
\begin{eqnarray*}
	y &=&\dfrac{y_{1}\left( r_{2}-r_{3}\right) -y_{2}\left( r_{1}-r_{3}\right)
		+y_{3}\left( r_{1}-r_{2}\right) }{x_{1}\left( r_{2}-r_{3}\right)
		-x_{2}\left( r_{1}-r_{3}\right) +x_{3}\left( r_{1}-r_{2}\right) }x \\
	&& \\
	&&+\dfrac{r_{1}\left( x_{2}y_{3}-x_{3}y_{2}\right) -r_{2}\left(
		x_{1}y_{3}-x_{3}y_{1}\right) +r_{3}\left( x_{1}y_{2}-x_{2}y_{1}\right) }{%
		x_{1}\left( r_{2}-r_{3}\right) -x_{2}\left( r_{1}-r_{3}\right) +x_{3}\left(
		r_{1}-r_{2}\right) }\text{.}
\end{eqnarray*}%
Similarly, we get the equation of the line passing through points $P_{13}$
and $P_{23}$ in the same way as the line above. Consequently, the points$\
P_{12}$, $P_{13}$ and $P_{23}\ $are collinear since the these points
determine one line called the axis of similitude (or Monge Line). Thus,
using the notations in \cite{Searby}, we can give the following main theorem.
\begin{thm}
	Let $C_{i}$\ and $C_{j}$ be two $\alpha -$circles (or $%
	l_{p}-$circles ) with centers \linebreak $M_{i}=\left( x_{i},y_{i}\right) $ and $%
	M_{j}=\left( x_{i},y_{i}\right) $ of radii $r_{i}$ and $r_{j}$ for $i,j\in
	\left\{ 1,2,3\right\} $, respectively. Then, the common tangents of the
	circles $C_{i}$\ and $C_{j}$ intersect at the point 
	\begin{equation*}
		P_{ij}=\left( \dfrac{r_{i}x_{j}-r_{j}x_{i}}{r_{i}-r_{j}},\dfrac{%
			r_{i}y_{j}-r_{j}y_{i}}{r_{i}-r_{j}}\right) \text{.}
	\end{equation*}%
	Also, the points $P_{12}$,$\ P_{13}$ and $P_{23}\ $of intersection of the
	three pairs of tangent lines lie on a line. This line is called Monge line
	whose equation%
	\begin{equation*}
		\left\vert 
		\begin{array}{ccc}
			x_{1} & x_{2} & x_{3} \\ 
			r_{1} & r_{2} & r_{3} \\ 
			1 & 1 & 1%
		\end{array}%
		\right\vert y=\left\vert 
		\begin{array}{ccc}
			y_{1} & y_{2} & y_{3} \\ 
			r_{1} & r_{2} & r_{3} \\ 
			1 & 1 & 1%
		\end{array}%
		\right\vert x-\left\vert 
		\begin{array}{ccc}
			x_{1} & x_{2} & x_{3} \\ 
			y_{1} & y_{2} & y_{3} \\ 
			r_{1} & r_{2} & r_{3}%
		\end{array}%
		\right\vert \text{.}
	\end{equation*}

\end{thm}




\begin{thebibliography}{1}
	\bibitem{Akca} Akca, Z., Kaya, R.: \textit{On the Distance Formulae In three Dimensional Taxicab Space.} Hadronic Journal. \textbf{27}, 521-532 (2006)
	
	\bibitem{Chen} Chen, G.: \textit{ Lines and Circles in Taxicab Geometry.} Master Thesis, Department of Mathematics and Computer Science, University of Central Missouri (1992)
	
	\bibitem{Colakoglu-Alpha} \c{C}olakoglu, H. B, Kaya, R.: \textit{A Generalization of Some Well-Known Distances and Related Isometries.} Math. Commun. \textbf{16}, 21-35 (2011)
	
	\bibitem{Thomas} Jahn, T., Spirova, M.: \textit{On bisectors in normed planes.} Contributions to Discrete Mathematics. \textbf{10}2, 1-9 (2015)
	
	
	
	
	
	\bibitem{Krause} Krause, E. F.: Taxicab Geometry. Addision-Wesley, Menlo Park, California (1975)
	
	\bibitem{Lassak} Lassak, M., Martini, H.: \textit{Reduced Convex Bodies in Finite Dimensional Normed Spaces.} Results. Math. \textbf{66}, 405–426 (2014)
	
	\bibitem{Leopol} Leopold, U., Martini, H.: Monge Points, Euler Lines, and Feuerbach Spheres in Minkowski Spaces, Geometry and Symmetry Conference,  Discrete Geometry and Symmetry, pp. 235-255 (2018)
	
	\bibitem{Swanepoel1} Martini, H., Swanepoel, K.J., Weiss, G.: \textit{The geometry of Minkowski spaces—a survey. Part I.} Expo. Math. \textbf{19}, 97–142 (2001)
	
	\bibitem{Swanepoel2} Martini, H., Swanepoel, K.J.: \textit{The geometry of Minkowski spaces—a survey. Part II.} Expo. Math. \textbf{22}, 93–144 (2004)
	
	\bibitem{Menger} Menger, K.: You Will Like Geometry. Guildbook of the Illinois Institute of Technology Geometry Exhibit, Museum of Science and Industry, Chicago (1952)
	
	\bibitem{Polar} Park, H. G., Kim, K. R., Ko, I. S., Kim,  B. H.: \textit{On Polar Taxicab Geometry In A Plane.} J. Appl. Math. Informatics. \textbf{32}, 783-790 (2014)
	
	\bibitem{Mustafaev} Mustafaev, Z., Martini, H.: \textit{On unit balls and isoperimetrices in normed spaces.} Colloquium Mathematicum. \textbf{127}, 133-142 (2012)
	
	\bibitem{Searby} Searby, D. G.: \textit{On Three Circles.} Forum Geometricorum. \textbf{9}, 181-193 (2009)
	
	\bibitem{Thompson} Thompson, A. C.:  Minkowski Geometry, Cambridge University Press, Cambridge (1996)
	
	\bibitem{Tian} Tian, S.:  \textit{Alpha Distance-A Generalization of Chinese Checker Distance and Taxicab Distance.} Missouri Journal of Mathematical Sciences. \textbf{17}(1), 35-40 (2005)
	
	\bibitem{Walker} Walker, W.:  \textit{Monge's Theorem in Many Dimensions.} Math. Gaz. \textbf{60}, 185-188  (1976)
	

\end{thebibliography}
\end{document}